\documentclass[11pt]{amsart}

\usepackage{amsmath}
\usepackage{amssymb}
\usepackage{amsthm}
\usepackage{amscd, amsfonts}
\usepackage[dvips]{epsfig}
\usepackage{verbatim}

\usepackage{epsf}
\usepackage[usenames]{color}
\usepackage[bookmarksnumbered,pdfpagelabels=true,plainpages=false,colorlinks=true,
            linkcolor=black,citecolor=black,urlcolor=blue]{hyperref}

\newcommand{\ddc}{d d^{c}}

\theoremstyle{plain}
\newtheorem{theorem}{Theorem}[section]

\theoremstyle{definition}

\numberwithin{equation}{section}

\title[A generalised Monge-Amp\`ere equation]{A generalised Monge-Amp\`ere equation}

\author[Pingali]{Vamsi P. Pingali}
\address{Department of Mathematics\\
Stony Brook University\\ Stony Brook, NY 11794, USA}
\email{Vamsi.Pingali@sunysb.edu}

\begin{document}
\maketitle

\begin{abstract}
We consider a generalised complex Monge-Amp\`ere equation on a compact K\"ahler manifold and treat it using the method of continuity. For complex surfaces,  we prove an easy existence result. We also prove that (for three-folds and a related real PDE in a ball in $\mathbb{R}^3$), as long as the Hessian is bounded below by a pre-determined constant (whilst moving along the method of continuity path), a smooth solution exists. Finally, we prove existence for another real PDE in a 3-ball, which is a local, real version of a conjecture of X.X.Chen. 
\end{abstract}

\section{Introduction}
Let $(X,\omega)$ be an $n$-dimensional, compact, K\"ahler manifold. Here, we consider a generalised complex Monge-Amp\`ere PDE (to be solved for a smooth function $\phi$)
\begin{gather}
\alpha_0(\omega + \ddc \phi)^n + \alpha _1 \wedge (\omega + \ddc \phi)^{n-1} + \ldots + \alpha _{n-1} \wedge (\omega + \ddc \phi) = \eta 
\label{maineq}
\end{gather}
where $\eta$, $\alpha _i$ are smooth, closed forms satisfying the obvious necessary condition $\int _X \eta = \int _X (\alpha _0 \omega ^n + \alpha _1 \wedge \omega ^{n-1} +\ldots)$. \\
\indent When $\eta >0$, $\alpha_i=0 \ \forall \ i \neq 0$ and $\alpha _0=1$, equation \ref{maineq} is the one introduced by Calabi \cite{Calabi} and solved by Aubin \cite{Aub} and Yau \cite{Yau}. Equations of this type are ubiquitous in geometry. A version of this generalised one appeared in \cite{XX}. The geometric applications of this equation shall be explored elsewhere. \\

\textbf{Acknowledgements}: The author thanks his adviser Leon A. Takhtajan for suggesting this direction of study and for sparing time generously to discuss the same. We also thank Xiu Xiong Chen, Dror Varolin, Yanir Rubinstein, and, Marcus Khuri for fruitful discussions.

\section{Statements of results}
We state a somewhat general theorem about uniqueness, openness and $C^0$ estimates. The proof is quite standard (adapted largely from \cite{Tian} which is in turn based on \cite{Yau}). Although the theorem is folklore, we haven't found the precise statement (in this level of generality) in the literature on the subject. In what follows, positivity of $(p,p)$ forms is strong positivity. Let $\mathcal{B}$ be the product of Banach submanifolds of forms wherein, an element of $\mathcal{B}$ is of the form $(\alpha _0, \alpha_1, \ldots, \alpha _{n-1}, \phi)$ where $\alpha_i$ are $C^{1,\beta}$ $(i,i)$, closed forms and $\phi$ is a $C^{3,\beta}$ function satisfying $\displaystyle \int _M \phi =0$, $n\alpha _0 (\omega + \ddc \phi)^{n-1} + (n-1) \alpha _1 \wedge (\omega + \ddc \phi)^{n-2} + \ldots + \alpha _{n-1}>0$ and, $\displaystyle \int _X (\sum _i \alpha _i \wedge \omega^{n-i}) \neq 0$. Also, let $\mathcal{\tilde{B}}$ be the Banach submanifold of $C^{1,\beta}$ top forms $\gamma$ with $\displaystyle \int _X \gamma =1$ and $\gamma >0$.
\begin{theorem}
If $\alpha _0 \omega ^n + \alpha _1 \omega ^{n-1} + \ldots >0$, $\eta >0$ and, $d\alpha _i =0$, then, any smooth solution $\phi$ of \ref{maineq} satisfying $\int _X \phi \omega ^n= 0$ and, $\kappa \geq K \omega ^{n-1}  $ where $K>0$ and $\displaystyle \sum _k (\alpha _k (\omega + \ddc \phi)^{n-k}-\alpha _ k\omega ^{n-k}) = \kappa \wedge \ddc \phi$, is bounded a priori: $\Vert \phi \Vert _{C^0} \leq C_{\eta}$. Also, if $\alpha _i >0 \ \forall i$ and, if there exists a smooth solution $\phi$ such that $\omega + \ddc \phi >0$, it is unique (upto a constant) among all such solutions; In addition, the mixed derivatives of $\phi$ are bounded a priori $\Vert \phi \Vert _{C^{1,1}} \leq C_{\eta}$. \\
\indent The map $T : \mathcal{B} \rightarrow \mathcal{\tilde{B}}$ defined by $T(\alpha _0, \alpha _1, \ldots, \phi) = \frac{\sum _i \alpha _i \wedge (\omega + \ddc \phi)^{n-i}}{\int _X (\sum _i \alpha _i \wedge \omega^{n-i})}$ is open and, so is the restriction of $T$ to a subspace defined by fixing the $\alpha_i$. Also, a level set of this map is locally a graph with $\phi$ being a function of the $\alpha _i$.  \\
\indent When $n=2$, and, $\alpha _0=1$, $\eta -\alpha _2 +\frac{\alpha _1 ^2}{4} >0$, there exists a unique, smooth solution to \ref{maineq} satisfying $\omega +\ddc \phi + \frac{\alpha _1}{2}>0$ .

\label{estimate}
\end{theorem}

In particular, if $\alpha _i = \omega ^i$ for some $i$ and all the other $\alpha _j$ are small enough, then, by the solution of the $k$-Hessian equations \cite{Hess}, \cite{Kol} we have a smooth solution of equation \ref{maineq}. \\ 
\indent One may formulate a version of the same problem locally as a Dirichlet problem on a pseudoconvex domain in $\mathbb{C}^n$. In this context, we note that viscosity solutions to the Dirichlet problem exist by \cite{LawsHarv} and \cite{Caff}. We also have the following result for a real version of the PDE :
\begin{theorem}
The following Dirichlet problem on the ball $B$ of radius $1$ centred at the origin 
\begin{eqnarray}
\det (D^2u) + \Delta u &=& tf + (1-t)36 \nonumber \\
u\vert _{\partial B}  &=& 0 \nonumber \\
f&>& 36 \nonumber
\end{eqnarray}
has a unique smooth solution at $t=T$ if $f$ is smooth and, for all $t\in [0, T)$, smooth solutions $u_t$ exist and satisfy $D^{2}u_t > 3$.  
\label{realv}
\end{theorem}
A similar result holds for complex three-folds.
\begin{theorem}
If $\alpha >0$, $\omega >0$ are smooth K\"ahler forms on a compact K\"ahler manifold $(X,\omega _0)$, then, there exists a constant $C>0$ depending only on $\alpha$ and $\omega _0$ such that, the equation
\begin{eqnarray}
(\omega + \ddc u_t)^3 + \alpha^2 (\omega + \ddc u_t) &=& \frac{e^{tf}\int (\omega ^3 + \alpha ^2 \omega)}{\int e^{tf} (\omega ^3 + \alpha ^2 \omega)}(\omega ^3 + \alpha ^2 \omega)  
\label{threeeq}
\end{eqnarray} 
has a unique smooth solution at $t=T$, if for all $t \in [0,T)$, smooth solutions exist and satisfy $\omega + \ddc u > C \omega _0$. 
\label{threefolds}
\end{theorem}
Finally, we present a local, real version of a conjecture of X.X. Chen (conjecture 4 in \cite{XX} made in the compact complex manifold case). Some progress has been made in a few special cases \cite{Lai}. However, in all these cases, the problem was reduced to an inverse Hessian equation. We prove existence in a special case here, using the method of continuity. Actually, a far more general result was proven in \cite{Kryl}, but, results on the Bellman equation were used (as opposed to a direct method of continuity). Such results may not carry over in an obvious way to the manifold case and hence our proof of this ``toy model''.  
\begin{theorem}
If $f>0$ is a smooth function on $\bar{B}(0,1)$ (the closed unit ball), then, the following Dirichlet problem has a unique, smooth, convex solution. 
\begin{eqnarray}
\det (D^2u) - \Delta u &=& f \nonumber \\
u\vert _{\partial B}  &=& 0 \nonumber \\
f&>& 0 \nonumber
\end{eqnarray}

\label{Chen}
\end{theorem}
\section{Standard results used}
For the convenience of the reader, we have included statements of some standard results in the form that we use in the proofs. \\
\indent Our principal tool to study fully nonlinear PDE like equation \ref{maineq}, is the method of continuity (It is like a flow technique. In fact this analogy was exploited more seriously, to great advantage, in  \cite{Yanir}). To solve $Lu =f$ where $L$ is a nonlinear operator, one considers the family of equations $Lu_t = \gamma (t)$ where, $\gamma (1)=f$ and $\gamma (0)=g$ such that, at $t=0$, one has a solution $Lu_0=g_0$. Then, one proves that the set of $t\in [0,1]$ for which the equation has a solution is both, open and closed (and clearly non-empty). In order to prove openness, one considers $L$ to be a map between appropriate Banach spaces. Then, the implicit function theorem of Banach spaces proves openness. However, while dealing with equations like Monge-Amp\`ere equations, one has to verify that certain conditions like ellipticity are preserved along the ``continuity path". This is crucial because, in order to solve the linearised equation and, to prove that indeed one has a solution in an appropriate Banach space, one needs ellipticity in these cases. In fact, in a few of the cases we shall consider, ellipticity is not preserved and hence, the best we can do is a ``short-time" existence result. In order to prove closedness, one needs to prove uniform (i.e. independent of $t$) \emph{a priori} estimates for $u$. In our case, we shall need these estimates in $C^{2, \alpha}$ in order to use the Arzela-Ascoli theorem to conclude closedness. These estimates are usually proved by improving on lower order estimates. Once, one produces a $C^{2,\alpha}$ solution, one ``bootstraps" the regularity (at teach $t \in [0,1]$) using the Schauder estimates. The Schauder estimates on a compact manifold (without boundary) are (they can be derived easily using similar interior and boundary ones in domains in $\mathbb{R}^n$ \cite{Trud}).
\begin{theorem}
\emph{Schauder $\emph{a priori}$ estimates on a Riemannian manifold}: If $Lu = f$, where $L$ is a second-order, uniformly elliptic operator with smooth coefficients, and, $u$ is a $C^{2,\alpha}$ and $f$ is a $C^{0,\alpha}$ function, 
\begin{gather}
\Vert u \Vert _{C^{2, \alpha}} \leq C(\Vert u \Vert _{C^{0}} + \Vert f \Vert _{C^{0, \alpha}}) \nonumber
\end{gather} 
\end{theorem} 
\indent In order to derive \emph{a priori} estimates, we shall use standard techniques as in \cite{Yau}, \cite{Tian} for the manifold case, and, \cite{Caff} for the Euclidean case. The main blackbox is the Evans-Krylov-Safanov theory for proving $C^{2, \alpha}$ estimates from $C^{2}$ ones. This requires (apart from uniform ellipticity) concavity of the equation. There is a similar version for the complex case. The real version is :
\begin{theorem}
Let $u$ be a smooth function on the unit ball satisfying, 
\begin{gather}
F(D^2 u, x, Du) = g \nonumber
\end{gather}
on the unit ball in $\mathbb{R}^n$ centred at the origin $B(0,1)$ with $u=0$ on the boundary of the ball. Here, $F$ is a smooth function defined on a convex open set of symmetric $n \times n$ matrices $\times \mathbb{R} \times \mathbb{R}^n$ which satisfies, \\
a) \emph{Uniform ellipticity on solutions} : There exist positive constants $\lambda$ and $\Lambda$ so that $0 < \lambda \vert \xi \vert^2 \leq F_{ij} (D^2 u ,x , Du) \xi _i \xi _j \leq \Lambda \vert \xi \vert ^2$ for all vectors $\xi$ and all $u$ satisfying the equation. \\
b) \emph{Concavity on a convex open set} : $F$ is a concave function on a convex open set of symmetric matrices (containing $D^2 u$ for all solutions $u$). \\
Then, $\Vert u \Vert _{C^{2, \alpha} (\bar{B(0,1)})} \leq C $ where $C$ and $\alpha$ depend on the first and second derivatives of $F$, $\Vert u \Vert _{C^2(\bar{B})}$, $\Vert g \Vert _{C^2(\bar{B})}$, $n, \lambda$ and $\Lambda$.    
\label{EvansKrylov}
\end{theorem}  
\indent The complex, interior version (that we need) is :
\begin{theorem}
Let $u$ be a $C^4$ function on the unit ball in $\mathbb{C}^n$ satisfying 
\begin{gather}
F(u _{i\bar{j}}, z,\bar{z}) = 0 \nonumber 
\end{gather}
for a $C^{2,\beta}$ function $F (x, p, \bar{p})$ satisfying,\\
a) \emph{Uniform ellipticity on solutions} : There exist positive constants $\lambda$ and $\Lambda$ so that $0 < \lambda \vert \xi \vert^2 \leq F_{i\bar{j}} (\ddc u, z, \bar{z}) \xi _i \xi _j \leq \Lambda \vert \xi \vert ^2$ for all vectors $\xi$ and all $u$ satisfying the equation. \\
b) \emph{Concavity on a convex open set} : $F$ is a concave function on a convex open set of hermitian matrices (containing $u_{i\bar{j}}$ for all solutions $u$). \\
Then, $\Vert u \Vert _{C^{2,\alpha}(B(0,\frac{1}{2}))} \leq C$ where, $C$ and $\alpha$ depend on $\lambda$, $\Lambda$, $n$, and $\Vert u_{i\bar{j}} \Vert _{C^0 (\bar{B})}$ and uniform bounds on the first and second derivatives of $F$ evaluated at $u$.
\label{CompEvans}
\end{theorem}
\indent The proofs are standard \cite{Trud}, \cite{Siu}, \cite{Blocki}, \cite{Kaz}. Usually, one proves these estimates when the PDE is concave on all symmetric matrices. In Monge-Amp\`ere equations, one needs a weaker requirement of being concave on a convex open set of symmetric matrices \cite{Caff} (the proofs go through easily with this requirement). \\
\indent To conclude, we add a few words about uniqueness. The usual technique for demonstrating uniqueness (due to Calabi) of $Lu=f$ is to assume two solutions $u_1$ and $u_2$, and, to write $0=Lu_1 - Lu_2 = \int _{0} ^{1} \frac{dL}{dt} (tu_2+(1-t)u_1) dt$. If the integrand is an elliptic operator, by the maximum principle, $u_1 = u_2$.  
\section{Proofs of the Theorems}

\subsection{Proof of theorem \ref{estimate}}
This proof is similar to the one for the usual Monge-Amp\`ere equation \cite{Tian}. \\
\emph{The $C^0$ estimate} : As usual, without loss of generality, we may change the normalisation to $\sup \phi = -1$ i.e. we may add $-1-\sup \phi$ to $\phi$. Indeed, if the new $\phi$ has a $C^0$ estimate, then, $\int _X \phi =0$ yields the desired $C^0$ estimate. This means, we just have to find a lower bound on $\phi$. Certainly $\phi$ has an $L^1$ bound \cite{Tian}. Let $\phi = -\phi _{-}$ (so that $\phi _{-} \geq 1$). Subtracting $\Theta = \displaystyle \sum_k \alpha _k \wedge \omega ^{n-k}$ and then, multiplying the equation by $\phi _{-} ^p$ and integrating, we have (here $\eta = e^f \Theta $),
\begin{eqnarray}
-\int \phi _{-}^p \ddc \phi_{-} \wedge \kappa &=& \int \phi _{-}^p (e^f-1)\Theta \nonumber \\
\int \phi _{-}^p (e^f-1)\Theta  &\leq& c \Vert \phi_{-} \Vert _{L^p} ^p \nonumber \\ 
-\int \phi _{-}^p \ddc \phi_{-} \wedge \kappa &=& \int d(\phi _{-} ^p) \wedge d^c \phi_{-} \wedge \kappa \nonumber \\
 &=& c\int d(\phi _{-} ^{\frac{p+1}{2}}) \wedge d^c (\phi _{-} ^{\frac{p+1}{2}}) \wedge \kappa \nonumber \\
&\geq& C\Vert \nabla (\phi _{-} ^{\frac{p+1}{2}}) \Vert _{L^2} ^2 \nonumber \\
&\geq& C_1 \left ( \int \phi _{-} ^{\frac{(p+1)n}{n-1}} \right )^{\frac{n-1}{n}} -C_2 \int \phi _{-} ^{p+1} \nonumber 
\end{eqnarray}
where the last inequality follows from the Sobolev embedding theorem. Upon rearranging, we have
\begin{eqnarray}
\Vert \phi _{-} \Vert _{L^{(p+1)(n)/(n-1)}} \leq (C (p+1))^{\frac{1}{p+1}} \Vert \phi _{-} \Vert _{L^{p+1}} \nonumber 
\end{eqnarray}
The Moser-iteration procedure gives $\sup \vert \phi \vert \leq C \Vert \phi \Vert _{L^2}$. If we prove that the right hand side is controlled by the $L^1$ norm of $\phi$, we will be done. Indeed,
\begin{eqnarray}
C\Vert \phi \Vert _{L^1} &\geq& \int \phi (1-e^f) \Theta  \nonumber \\
&\geq& C \Vert \nabla \phi \Vert _{L^2} ^2 \nonumber \\
&\geq& C \Vert \phi - \langle \phi \rangle \Vert _{L^2} ^2 \nonumber \\
\Vert \phi \Vert _{L^2} &\leq& C(\Vert \phi \Vert _{L^1} +1) \nonumber  
\end{eqnarray}
where, we have used the Poincar\' e inequality. 
Hence, proved. \\
\newline
\emph{Uniqueness} : If $\phi _1$ and $\phi _2$ are two solutions, upon subtraction we have, 
\begin{gather}
\sum _i \alpha_i \wedge ((\omega + \ddc \phi _1)^{n-i} - (\omega + \ddc \phi_2)^{n-i}) = 0 \nonumber \\
\Rightarrow  \int _{0}^{1} \displaystyle \sum _ k k\alpha_k  \wedge (\omega + \ddc \phi_1+ t\ddc (\phi _2-\phi_1))^{n-k-1}dt  \wedge \ddc (\phi_2 - \phi_1)  = 0 \nonumber 
\end{gather}
 Thus, by the maximum principle, $\phi_2 - \phi_1$ is a constant. \\
\newline
\emph{The mixed derivatives estimate}: When $\alpha _i >0$, 
\begin{gather}
\eta \geq \alpha _{n-1} \wedge (\omega + \ddc \phi)
> C(\mathrm{tr}(\omega + \ddc \phi))   \nonumber 
\end{gather}
where $C>0$. Since $0<\omega + \ddc \phi$, the eigenvalues of $\ddc \phi$ are bounded above. Thus, the mixed second derivatives of $\phi$ are bounded. Note that, by the Schauder estimate \cite{Morr}, the first derivatives are bounded as well.  \\
\newline
\emph{Openness} :  The map $T$ is smooth. Its G\^ateaux derivative is $DT(0, 0, \ldots, 0, \chi) = (n\alpha _0 (\omega + \ddc \phi)^{n-1} + (n-1) \alpha _1 \wedge (\omega + \ddc \phi)^{n-2} + \ldots + \alpha _{n-1}) \wedge \ddc \chi$. It is clearly a bounded surjection (by the Schauder theory) onto its image if $n\alpha _0 (\omega + \ddc \phi)^{n-1} + (n-1) \alpha _1 \wedge (\omega + \ddc \phi)^{n-2} + \ldots + \alpha _{n-1} >0$. If $DT$ is restricted to vectors of the form $(0,0,\ldots,0,\chi)$, then, it is a Banach space isomorphism. Hence, by the implicit function theorem of Banach manifolds, openness is guaranteed. In fact, it also guarantees that, on a level set, $\phi$ can be solved for (locally), in terms of $\alpha_i$.  \\
\newline
\emph{The n=2 case} : The equation we have is equivalent to
\begin{gather}
(\omega + \ddc \phi + \frac{\alpha_1}{2}) ^2 = \eta - \alpha _2 + \frac{\alpha _1 ^2}{4} \nonumber  
\end{gather}
This is just the usual Monge-Amp\`ere equation and hence we are done.
\subsection{Proof of theorem \ref{realv}}
Uniqueness is proved as before. We shall only prove existence. Let $Lu = \det (D^2 u) + \Delta u$. To this end, we use the method of continuity. Consider the equation
\begin{eqnarray}
Lu_t &=& tf +(1-t)L\phi \nonumber \\
u_t \vert _{\partial B} &=& 0 \nonumber \\
\phi &=& \frac{3}{2} \sum x_i^2 -\frac{3}{2} 
\label{tequa}
\end{eqnarray}
When $t=0$, it has a smooth solution, namely, $\phi$. \\
\emph{Openness}: Let $\Omega \subset C_{0}^{2,\alpha} (\bar{B})$ be the set of $u$ such that $D^2 u > 3$ (where the subscript $0$ indicates vanishing on the boundary). This is an open subset. Define $T: \Omega \rightarrow C_{0} ^{0,\alpha}$ to be $T(u_t) = \det (D^2 u_t) + \Delta u_t$. If $u_s$ is a solution of \ref{tequa}, then, it is easy to see that $DT_{u_s}$ is a linear isomorphism. Hence, by the inverse function theorem of Banach manifolds, we see that the set of $t$ for which there is a solution is open. \\ 
\emph{Closedness}: Suppose there is a sequence $t_i \rightarrow t$ such that there are smooth solutions $u_{t_i}$ satisfying $D^2 u > 3$. Then, we wish to prove that a subsequence of the $u_{t_i}$ converges to a smooth solution $u_t$ in the $C^{2,\beta}$ topology. This requires \emph{apriori} estimates (the convergence following from the Arzela-Ascoli theorem). We shall prove the same for the equation \ref{realv}. We just have to prove the $C^{2,\alpha}$ estimate in order to ensure smoothness (by the Schauder theory). \\
\emph{$C^0$ estimate}: Note that $\Delta u \leq f$. Hence, for $A >>1$, $0> f-\Delta(A \sum x_i^2) = \Delta (u-A\sum x_i ^2)$. The minimum principle implies that $u \geq A\sum x_i ^2- A$. Since, $\Delta u > 9$, $u\leq 0$ by the maximum principle. thus, $\Vert u \Vert _{C^0} \leq C$.  \\
\emph{$C^1$ estimate}: Differentiating both sides using the operator $D$, $\mathrm{tr}((\mathrm{Hess} u)^{-1} D^2w) + \Delta w = Df$ where $w=Du$. Just as before, by adding or subtracting a large multiple of $\sum x_i ^2$ to $w$ and using the maximum principle, we see that $\Vert Du \Vert _{C^0}$ is controlled by its supremum on the boundary. The tangential boundary derivatives are $0$. Since, $A \sum x_i ^2 - A \leq u \leq 0$, $ \vert \frac{\partial u}{\partial n} \vert \leq 2A $. Hence, $\Vert u \Vert _{C^1} \leq C$.  \\
\emph{$C^2$ estimate}:  Since $\Delta u \leq Lu \leq f$ and $\Delta u>0$, $\Vert u_{ij} \Vert _{C^0} \leq C$. Hence, $\Vert u \Vert _{C^2} \leq C$.  \\
\emph{$C^{2,\alpha}$ estimate}: So far, we haven't used anything about the sequence except that $D^2 u_{t_i} >0$. This will change presently. For any function $F: \mathbb{R} \rightarrow \mathbb{R}$, $F(\det (D^2 u) + \Delta u) = F(f)$. If we choose the function appropriately, then the resulting equation will be a concave, uniformly elliptic Monge Amp\`ere PDE to which we may apply the Evans-Krylov theory to extract a $C^{2, \alpha}$ estimate.  \\
\indent We claim that, the function $F(x) = \int _{36} ^{x} e^{-\frac{t^2}{2}} dt$ is such that, $g(\lambda _1, \lambda _2, \lambda _3) = F(\sum \lambda _i + \lambda _1 \lambda _2 \lambda _3)$ has a uniformly positive gradient and is concave if $\lambda _i > 3$.  By using the $C^2$ estimate and theorem \ref{EvansKrylov}, we have the desired estimate. \\
\indent We shall prove the aforementioned fact:  Let $x=\sum \lambda _i + \lambda _1 \lambda _2 \lambda _3$. We see that $\frac{\partial g}{\partial \lambda _i}\vert_{D^2 u} = e^{-\frac{x^2}{2}}(1+\frac{\lambda _1 \lambda _2 \lambda _3}{\lambda _i}) > e^{-\frac{f^2}{2}}$ and is less than $1+3f$ where we have evaluated the derivative at the eigenvalues of the Hessian of a solution of equation \ref{realv}. Hence, it is uniformly elliptic. \\
If $(v_1 , v_2 , v_3 ) \in \mathbb{R}^3 $, then $ -v_i v_j\frac{\partial ^2 g}{\partial x_i \partial x_j}  =  e^{-\frac{x^2}{2}}(x(v_1(1+\lambda _2 \lambda _3)+v_2 (1+\lambda _3 \lambda _1) + v_3 (1+\lambda _1 \lambda _2))^2)-2e^{-\frac{x^2}{2}}(v_1 v_2 \lambda_3 + v_2 v_3 \lambda _1 + v_3 v_1 \lambda _2)$, which is in turn equal to 
$ e^{-\frac{x^2}{2}} (v_1 ^2 \alpha + \beta v_1 + \gamma)  \geq 0 \Leftrightarrow \beta ^2 - 4 \alpha \gamma \leq 0$ and, 
\begin{eqnarray}
\frac{\beta ^2- 4\alpha \gamma}{4} &=&  (v_2(x(1+\lambda _2 \lambda _3)(1+\lambda _1 \lambda _3)-\lambda _3)+v_3 (x(1+\lambda _1 \lambda _2)(1+\lambda _2 \lambda _3)-\lambda _2))^2 \nonumber \\
&-& x(1+\lambda _2 \lambda _3)^2(v_2 ^2 x (1+\lambda _1 \lambda _3)^2 + v_3 ^2 x(1+\lambda _1 \lambda _2)^2 
 + 2v_2 v_3 (x(1+\lambda _1 \lambda _3)(1+\lambda _1 \lambda _2)-\lambda _1)) \nonumber \\
&=& \tilde{\alpha}v_2 ^2 + \tilde{\beta} v_2 + \tilde{\gamma}  \leq 0 \nonumber  
\end{eqnarray}
with the last inequality holding if and only if $\tilde{\alpha} \leq 0$ and $\tilde{\beta}^2 - 4\tilde{\alpha} \tilde{\gamma}\leq 0$. Let us assume (without loss of generality) that $v_3 \neq 0$ and that $\lambda _1 < \lambda _2 <\lambda _3$. 
\begin{eqnarray}
\tilde{\alpha} &=& (x(1+\lambda _2 \lambda _3)(1+\lambda _1 \lambda _3)-\lambda _3)^2-x^2(1+\lambda _2 \lambda _3)^2(1+\lambda _1 \lambda _3)^2 \nonumber \\
&=& \lambda _3 ^2 - 2\lambda _3 x (1+\lambda _2 \lambda _3)(1+\lambda _1 \lambda _3) \nonumber \\
&\leq& -2x \lambda _3 ^2 (\lambda _1 +\lambda _2 + \lambda _1 \lambda _2 \lambda _3)
\leq -2x \lambda _3 ^2 \frac{2x}{3}\nonumber 
\end{eqnarray}
\begin{eqnarray}
\frac{\tilde{\gamma}}{v_3 ^2} &=&(x(1+\lambda _2 \lambda _1)(1+\lambda _2 \lambda _3)-\lambda _2)^2-x^2(1+\lambda _2 \lambda _3)^2(1+\lambda _1 \lambda _2)^2 \nonumber \\
&\leq& -\frac{4x^2\lambda _2^2}{3}\nonumber 
\end{eqnarray}
\begin{eqnarray}
\frac{\tilde{\beta}}{2v_3} &=& (x(1+\lambda _2 \lambda _3)(1+\lambda _1 \lambda _2)-\lambda _2)(x(1+\lambda _2 \lambda _3)(1+\lambda _1 \lambda _3)-\lambda _3)\nonumber \\
&-&x(1+\lambda _2 \lambda _3)^2(x(1+\lambda _1 \lambda _3)(1+\lambda _1 \lambda _2)-\lambda _1)\nonumber 
\end{eqnarray}
\begin{eqnarray}
\left (\frac{\tilde{\beta}}{2v_3} \right ) ^2 &=& (x(1+\lambda _2 \lambda _3)(\lambda _2 + \lambda _3 - \lambda _1 + \lambda _1 \lambda _2 \lambda _3)-\lambda _2 \lambda _3)^2 \nonumber \\
&\leq& x^4 (1+\lambda _2 \lambda _3) ^2 \nonumber \\
\frac{\tilde{\beta} ^2 - 4 \tilde{\alpha}\tilde{\gamma}}{4v_3 ^2} &\leq& 0 \nonumber 
\end{eqnarray}

Hence proved. \\
\emph{Remark} : Writing equation \ref{maineq} for $n=3$ and $\alpha_0 =1$ we have, 
\begin{gather}
(\omega + \ddc \phi)^3 + \alpha _1 (\omega + \ddc \phi) ^2 + \alpha _2 (\omega + \ddc\phi) = \eta \nonumber \\
\Rightarrow (\omega + \frac{\alpha _1}{3} + \ddc \phi)^3 + (\alpha _2 - \frac{\alpha _1 ^2}{3})(\omega + \frac{\alpha _1}{3} + \ddc \phi) = \eta - \frac{2\alpha _1 ^3}{27} + \frac{\alpha _1 \alpha _2}{3} \nonumber 
\end{gather}
A local, real version of a special case of the above is equation \ref{realv}.
\subsection{Proof of theorem \ref{threefolds}}
 Once again, we apply the method of continuity. We shall impose several conditions on $C$ (as we go along). It should be large enough so that, whenever $\beta > C\omega _0$, $\beta ^3 > 3\alpha ^2 \beta$ (Indeed, if $K>0$ and $B>0$ are given, $\det(A) > K\mathrm{tr}(BA)$ for sufficiently large $A>0$). Obviously, at $t=0$, $u=0$ solves the equation. Openness and uniqueness, follow from theorem \ref{estimate}. As before, if $t_i \rightarrow t$ is a sequence such that there exist smooth solutions $u_i$ satisfying $\omega + \ddc u_i > C \omega _0$, then, we shall prove that a subsequence converges to a smooth solution $u$ in the $C^{2,\beta}$ topology. As usual, we need \emph{apriori} estimates for this.\\
\indent The $C^{0}$ and the mixed derivative estimates follow directly from theorem \ref{estimate}. We have to prove the $C^{2,\alpha}$ estimate (thus proving existence and smoothness as before). It suffices to prove a local (interior) estimate. We shall accomplish this via the complex version of the (interior) Evans-Krylov theory done in \cite{Blocki} and \cite{Siu}. \\
\indent The local (in a ball) version of the equation is
\begin{eqnarray}
\det (\phi _{i \bar{j}}) + \mathrm{tr}(B^{-1}[\phi _{i\bar{j}}]) &=& f \nonumber \\
\phi _{i \bar{j}} &>& C>1 \nonumber \\ 
f &>& C^3+9\Vert B ^{-1} \Vert^2 C   \nonumber \\
 \label{locver} 
\end{eqnarray} 
where $B^{-1}_{i\bar{j}} = \det(\alpha) [\alpha]^{-1} _{i\bar{j}}$. We claim that the function $g (A) = F(\det (A) + \mathrm{tr}(B^{-1}A))$ from hermitian matrices satisfying $A > CId $ to $\mathbb{R}$ (where $F(x) = \int _{c} ^{x} e^{-\frac{t^2}{2}}dt$) is concave and uniformly elliptic. Let the eigenvalues of $A$ be $\lambda _1, \lambda _2$ and $\lambda _3$. The uniform ellipticity is trivial (as in the proof of theorem \ref{realv}). The concavity is also somewhat similar to theorem \ref{realv}, but requires some modification. Indeed (here $V$ is an arbitrary hermitian matrix and $x=\det (A) + \mathrm{tr}(B^{-1}A)$),
\begin{eqnarray}
g^{''}(V,V) &=& g''(x) (\det (A) \mathrm{tr}(A^{-1}V)+\mathrm{tr}(B^{-1}V))^2 \nonumber \\
 &+& g'(x) (-\det(A)\mathrm{tr}((A^{-1}V)^2)+\det(A)(\mathrm{tr}(A^{-1}V))^2) \nonumber 
\end{eqnarray}
\indent We wish to prove that $g^{''}(V,V)<0$ for every hermitian $V$. Let's diagonalise the positive-definite form $B^{-1}$, i.e. $PB^{-1}P^{\dag} = I$ for some matrix $P$. Define $\tilde{A} = (P^{\dag})^{-1} A P^{-1}$ and $\tilde{V} = (P^{\dag})^{-1} V P^{-1}$. Now, using a unitary matrix $U$, we may diagonalise $\tilde{A}$ i.e. $\tilde{\tilde{A}} = U \tilde{A} U^{\dag} = \mathrm{diag} (a_1, a_2, a_3)$ where $a_1 \leq a_2 \leq a_3$ and $\tilde{\tilde{V}} = U \tilde{V} U^{\dag}$. This implies that $\det (\tilde{A}) \det(B) = \det(A)$ and $\mathrm{tr} ( \tilde{A}) = \mathrm{tr}(B^{-1}A)$. Let $\tilde{\tilde{V}}_{ii} = v_i$.  Hence,
\begin{eqnarray}
g^{''}(V,V) &=& -xe^{-\frac{x^2}{2}} (\det(B)a_1 a_2 a_3 (\sum \frac{ v_i}{a_i}) + \sum v_i) ^2 \nonumber \\
&+& e^{-\frac{x^2}{2}} ((\sum \frac{ v_i}{a_i})^2-(\sum \frac{v_i ^2}{a_i ^2} + 2 (\frac{\vert v_{12} \vert^2}{a_1 a_2}+\frac{\vert v_{23} \vert^2}{a_2 a_3}+\frac{\vert v_{13} \vert^2}{a_1 a_3}))) \nonumber \\
&\leq&  -xe^{-\frac{x^2}{2}} (\sum v_i (\det(B)\frac{a_1 a_2 a_3}{a_i}  + 1)) ^2 
+ 2\det(B)e^{-\frac{x^2}{2}} (v_1 v_2 a_3 + v_2 v_3 a_1 + v_3 v_1 a_2) \nonumber \\
&=& e^{-\frac{x^2}{2}}(Pv_1 ^2 + Q v_1 +R) \nonumber 
\end{eqnarray}
where,
\begin{eqnarray}
P &=& -x (\det(B) a_2 a_3 +1)^2 \leq 0\nonumber \\ 
Q &=& 2(\det(B)(v_2 a_3 + v_3 a_2)-x (\det(B)a_2 a_3 + 1)(v_2 (\det(B) a_1 a_3 +1) + v_3 (\det(B) a_1 a_2 + 1))))\nonumber \\
R &=& 2 \det(B) v_2 v_3 a_1 - x (v_2 (\det(B) a_1 a_3 + 1) + v_3 (\det(B) a_1 a_2 + 1))^2\nonumber 
\end{eqnarray}
as before, we want $Q^2 - 4 PR < 0$. Assume (without loss of generality) that $v_3=1$.  
\begin{eqnarray}
\frac{Q^2-4PR}{4} &=& Jv_2 ^2 + K v_2 + L \nonumber 
\end{eqnarray} 
where,
\begin{eqnarray}
J &=& \det(B)  a_3  (\det(B)a_3 - 2x(\det(B)a_2 a_3 + 1)(\det(B) a_1 a_3 +1))\nonumber \\
&<& \det(B) ^2 a_3 (1-2a_1 a_2 \det (B)) \nonumber \\
&<& 0 \nonumber \\
K &=& 2 (a_2 a_3 (\det(B))^2 -x(\det(B)a_2 a_3 +1)((\det(B)) ^2 a_1 a_2 a_3 + \det(B) (a_2+a_3-a_1)))\nonumber \\
L &=& \det(B)a_2 (\det(B)a_2 - 2x(\det(B)a_2 a_3 + 1)(\det(B) a_1 a_2 +1)) \nonumber 
\end{eqnarray}
where the first inequality follows from the assumption that $\det(A) =\det(B) a_1 a_2 a_3 > 3\mathrm{tr}(B^{-1}A) > 3\sum a_i$ . Hence,
\begin{eqnarray}
\frac{K^2}{4} &=& (a_2 a_3 (\det(B))^2 -x\det(B)(\det(B)a_2 a_3 +1)(\det(B) a_1 a_2 a_3 + (a_2+a_3-a_1)))^2\nonumber \\
&\leq& x^4(\det(B) a_2 a_3 + 1)^2 \det(B) ^2 \nonumber \\
J &\leq& -2x \det(B) ^2 a_3 ^2 (\det(B) a_1 a_2 a_3 + a_1 + a_2) \nonumber \\
&\leq& -2x \det(B) ^2 a_3 ^2 \frac{2x}{3} \nonumber \\
K &\leq& -2x \det(B) ^2 a_2 ^2\frac{2x}{3} \nonumber  
\end{eqnarray}
Thus, $K^2-4JL < 0$ implying that $g$ is concave. The $C^{2, \alpha}$ estimate follows from theorem \ref{CompEvans}. 

\subsection{Proof of theorem \ref{Chen}}
We use the method of continuity again. As before, openness follows easily using the Implicit function theorem on Banach spaces. Here, we prove only the \emph{a priori} estimates. Smoothness follows by bootstrapping, as indicated earlier. Lastly, we shall also prove the uniqueness of convex solutions. \\
\emph{$C^0$ estimate}: Since $u$ is convex, its maximum is attained on the boundary and hence, $u\leq 0$. Let $\phi = \frac{\mu}{2} r^2 - \frac{\mu}{2}$ where $\mu >0$, and, $\mu^3-3\mu > \max f$; then, subtracting $\det(D^2 u) - \Delta u$ from $\det(D^2 \phi) - \Delta \phi$, we have (assume that the eigenvalues of $D^2 u$ are $\lambda _i$),
\begin{gather}
L(\phi - u) = \det(D^2 \phi) - \det (D^2 u) - \Delta (\phi - u) \nonumber   \\
 = (\mu - \lambda _1)(\frac{\mu ^2 + \frac{\mu}{2} (\lambda _2 + \lambda _3)+\lambda _2 \lambda _3}{3}-1) + (\mu - \lambda _2)\ldots \nonumber \\
= \mu ^3 - \mu - f > 0 \nonumber 
\end{gather}    
\indent We see that, since $\mu^2 >3$, hence, $L$ is an elliptic operator acting on $\phi -u$ with $L (\phi - u)>0$. So, by the maximum principle, $\phi < u$. This gives us a $C^0$ estimate on $u$.\\
\emph{$C^1$ estimate}: Follows from ellipticity as before.\\
\emph{$C^2$ estimate}: For future use, notice that, atleast two of the eigenvalues of $D^2 u$ are larger than $1$. Taking derivatives of the equation we have (let $u_0$ be the minimum of $u$), 
\begin{gather}
\det(D^2 u) \mathrm{tr}((D^2 u)^{-1} D^2 u_i) - \Delta u_i = f_i \nonumber  \\
\det(D^2 u) \mathrm{tr}((D^2 u)^{-1} D^2 \Delta u) - \Delta \Delta u = \Delta f + \sum _i \det (D^2 u) \mathrm{tr}(((D^2 u)^{-1} D^2 u_i)^2 ) \nonumber \\
 - \det (D^2 u) \sum _i (\mathrm{tr}((D^2 u)^{-1} D^2 u_i))^2  \nonumber   
\end{gather}
\indent Let $A = \det (D^2 u) (D^2 u)^{-1} - I$. Consider $g = \Delta u + \mu (u-u_0)>0$ (we shall choose the constant $\mu >0$ later. It can depend on $\Vert u \Vert _{C^1} $ and other constants). Notice that, if $g$ is bounded, then, so is $\Delta u$ and thus, $D^2 u$ is bounded. At the maximum of $g$ (if it occurs in the interior), $(\Delta u)_i = -\mu u_i$ and $\mathrm{tr}(A D^2 g) \leq 0$. This implies,
\begin{eqnarray}
0 &\geq& \Delta f + \mu\mathrm{tr} (A D^2 u) - \frac{\vert -\mu \nabla u + \nabla f \vert^2}{\Delta u + f} \nonumber \\
&\geq& C_1 (\mu) - \frac{C_2 (\mu)}{\Delta u + f} + (2\Delta u +3f)\mu \nonumber 
\end{eqnarray}  
\indent Hence, $\Delta u$ is bounded at that point. Thus, $g$ is bounded at that point. This implies that $\Delta u$ is bounded everywhere. If the maximum of $g$ occurs on the boundary (call the max $g_0$), we shall have to analyse it separately. Let $\tilde{g} = g + g_0 (1-2r^2)$. Clearly, the maximum of $\tilde{g}$ has to occur in the interior. There, $D\tilde{g}=0$ and $\mathrm{tr}(AD^2\tilde{g}) \leq 0$. Hence (here, we assume that, $\mathrm{tr}(A) = \sum (\lambda _i \lambda _j - 1) $ and, that $g_0$ are sufficiently large compared to constants; If not, we are done), 
\begin{eqnarray}
0 &\geq& -C_1 (\mu) + \mu \mathrm{tr}(A) g_0 -C_2 g_0 - (C_4 g_0 + C_6 (\mu)) \mathrm{tr}(A) - C_5 \frac{g_0 ^2}{\Delta u +f} \nonumber \\
&\geq& -\tilde{C}_1 (\mu) + \mu \mathrm{tr}(A) g_0 -\tilde{C}_2 g_0 - \tilde{C}_4 g_0 \mathrm{tr}(A) - \tilde{C}_5 g_0 \nonumber \\
&\geq& -E_1 (\mu) + (\mu g_0 - E g_0) \mathrm{tr}(A) \nonumber 
\end{eqnarray} 
Choosing $\mu >E$, we see that, $g_0$ is bounded. Notice that, this also implies a lower bound on $D^2 u$. This is because, $M\lambda_i > \lambda _1 \lambda _2 \lambda _3 > f$. \\
\emph{$C^{2, \alpha}$ estimate} : Notice that, the set $Y$ of positive, symmetric matrices satisfying $\det(A)-\mathrm{tr}(A)>0$ is a convex open set (lemma $4.16$ of \cite{Kryl}). Also, our equation maybe written as $-1=-\frac{\Delta u}{\det (D^2 u)} - \frac{f}{\det(D^2 u)} = F(D^2 u, x)$ which is certainly concave on $Y$ by the same lemma in \cite{Kryl}. It is uniformly elliptic on solutions as long as the eigenvalues of the Hessian are bounded below and above (which they are, by the $C^2$ estimates). Theorem \ref{EvansKrylov} yields the desired estimates.\\
\emph{Uniqueness} :  If $u_1$ and $u_2$ are two convex solutions of the equation $F(u)=-1$ (as above), then, upon subtraction, $ 0 = \int _{0} ^{1} \mathrm{tr}((-I+\det (D^2 u_t) (D^2 u_t)^{-1})D^2 (u_2-u_1)) dt = L(u_2-u_1)$ where, $L$ is elliptic. By the maximum principle, $u_1 = u_2$.

\end{document}